\documentclass[a4paper]{amsart}
\usepackage{amsthm,mathtools,amssymb}
\usepackage{graphicx}
\usepackage{tikz-cd}
\usepackage{mathrsfs}
\usepackage[colorlinks=false, linkcolor=red]{hyperref}
\usepackage{enumerate}
\usepackage{cleveref}

\theoremstyle{plain}
\newtheorem{theorem}{Theorem}[section]
\newtheorem{conjecture}[theorem]{Conjecture}

\newtheorem{proposition}[theorem]{Proposition}
\newtheorem{corollary}[theorem]{Corollary}
\newtheorem{lemma}[theorem]{Lemma}
\newtheorem{remark}[theorem]{Remark}

\newcommand{\et}{\mathrm{et}}

\title[]{Cylinder maps of algebraic cycles on cubic hypersurfaces}

\author {Renjie Lyu}
\address{School of Mathematical Sciences, Xiamen University, Xiamen 361005, China.}
\email{r.lyu@xmu.edu.cn}
\date{\today}
\thanks{This work is adapted from the author's doctoral dissertation, ``Algebraic cycles on cubic hypersurfaces and Fano scheme of lines,'' Universiteit van Amsterdam, 2020. The program was supported by De Nederlandse Organisatie voor Wetenschappelijk Onderzoek, 614.001.507.}

\begin{document}
\begin{abstract}
Let \(X\subset \mathbb{P}^{n+1}\) be a smooth cubic hypersurface, and 
let \(F(X)\) be the variety of lines on \(X\). We prove the surjectivity 
of the cylinder maps on the Chow groups of \(X\) and \(F(X)\) if \(X\)
contains a one-cycle of degree \(1\). Mongardi and Ottem previously
proved the integral Hodge conjecture for curve classes on hyperk\"ahler 
manifolds. Using the cylinder maps, we provide an alternative proof 
for the \(F(X)\) of a smooth complex cubic fourfold \(X\), which is a 
special hyperk\"ahler fourfold. In addition, we confirm the integral Tate 
conjecture for \(F(X)\) of a smooth cubic fourfold \(X\) over a finitely 
generated field.
\end{abstract}

\maketitle
\tableofcontents

\section*{Introduction}
Let \(X\) be a smooth closed subvariety in the projective space 
\(\mathbb{P}^N_k\) over a field \(k\). The \emph{variety of lines} 
on \(X\), denoted by \(F(X)\), is the parameter space of lines in 
\(\mathbb{P}^N_k\) that are contained in \(X\). The geometry and topology
of \(F(X)\) and \(X\) are closely related, particularly when \(X\) is 
a Fano variety.

Fix an integer \(r\geq 1\). The \emph{cylinder map} or 
\emph{cylinder homomorphism} associated to lines on \(X\) is a
homomorphism on groups of algebraic cycles:
  \[
  \psi\colon \mathcal{Z}_{r-1}(F(X))\to \mathcal{Z}_{r}(X)
  \]
such that for a given \((r-1)\)-dimensional cycle \(\gamma\) on \(F(X)\),
the \(r\)-dimensional cycle \(\psi(\gamma)\) on \(X\) is the union of 
lines over \(\gamma\). On Chow groups, this map is naturally described
via the universal \(\mathbb{P}^1\)-bundle over \(F(X)\), namely, 
the incidence correspondence     
  \[
  P\coloneqq \{([\ell], x)\in F(X)\times X~|~ x\in \ell\subset X\}.
  \]   
The cylinder map on the Chow groups is then given by the induced homomorphism
  \begin{equation}
    \label{intro:cylinder-homo}
  [P]_*=q_*p^*\colon \mathrm{CH}_{r-1}(F(X))\to \mathrm{CH}_{r}(X). 
  \end{equation}
where \(p\colon P\rightarrow F(X)\) and \(q\colon P\rightarrow X\)
are the natural projections.

Suppose that \(k=\mathbb{C}\). The cohomological cylinder map
  \begin{equation}
  \label{eqs:coh-cyl-map}
  q_*p^*\colon \mathrm{H}_{n-2}(F(X), \mathbb{Z})\rightarrow \mathrm{H}_{n}(X, \mathbb{Z})
  \end{equation}
where \(n=\dim X\) has been studied in various cases, especially 
for Fano threefolds~\cite{Welter:Abel-Jacobi} and related variants~\cite{Letizia:Abel-Jacobi}. 
When \(X\) is a smooth complex cubic threefold, Clemens and Griffiths
~\cite{Clemens-Griffith:cubic3fold} proved that~\eqref{eqs:coh-cyl-map} 
is an isomorphism, and this fact played a crucial role in the proof of
Torelli and irrationality theorems for cubic threefolds. For a smooth
complex cubic fourfold \(X\), Beauville and Donagi~\cite{Beauville-Donagi:cubic4fold} 
showed that \(F(X)\) is a hyperk\"ahler \(4\)-fold and~\eqref{eqs:coh-cyl-map} 
is an isomorphism, which is fundamental in describing the 
Bogomolov-Beauville-Fujiki form on \(\mathrm{H}^2(F(X), \mathbb{Z})\). 
More generally, Shimada asserted that~\eqref{eqs:coh-cyl-map} is surjective 
for general Fano complete intersections~\cite{Shimada:cylin-hom-Fano-ci},
and Voisin generalized the conclusion to any smooth Fano complete 
intersections~\cite{Voisin:coniveauRC}.

In this paper, we concern the cylinder maps for a smooth cubic 
hypersurface. Our main result is as follows:
\begin{theorem}[=Theorem~\ref{surj-cylinder}]
  \label{intro:main-thm-surj}
Let \(k\) be a field, and let \(X\subset \mathbb{P}^{n+1}_k\) be a smooth
cubic hypersurface of dimension \(n\geq 3\). Let \(F(X)\) be the
variety of lines on \(X\). Assume that \(X\) contains a one-cycle of 
degree \(1\). Then for any integer \(r\geq 1\), the cylinder map
\[
q_*p^*\colon \mathrm{CH}_{r-1}(F(X))\rightarrow \mathrm{CH}_{r}(X)
\]
is surjective except \(r=\dim X-1\).
\end{theorem}
It was a question whether any one-cycle on a smooth cubic hypersurface
\(X\) is rationally equivalent a finite sum of lines, which 
is equivalent to the surjectivity of the cylinder map for \(r=1\). 
If \(k\) is an algebraically closed field, Paranjape~\cite{Paranjape94} 
answered this question affirmatively for \(\dim X\geq 5\). 
Tian and Zong studied one-cycles on separably rationally connected 
variety~\cite{Tian-Zong}. In particular, they proved that the Chow group
of one-cycles on a separably rationally connected Fano complete 
intersection of index at least \(2\) is generated by lines~\cite[Thm. 1.7]{Tian-Zong}. 

Very few examples exist where the cylinder maps of higher dimensional 
cycles are known to be surjective. Mboro proved the surjectivity for 
\(r\)-cycles on a smooth cubic hypersurface with \(r\geq 1\) and 
\(n\geq 2^r+1\), under the assumption that resolution of singularities 
exists in dimension \(\leq r\), see~\cite[Thm. 0.2]{Mboro:CH2}.
Additional examples can be found in Lewis~\cite{Lewis:cyl-map}. 

Our proof of Theorem~\ref{intro:main-thm-surj} follows Shen's method 
in~\cite{shen14,Shen-uni-gen} where he estabished key relations among 
one-cycles on a cubic hypersurface over an arbitrary field. These 
relations can deduce the surjectivity of the cylinder map for \(r=1\). 
For higher dimensional cycles, we obtain similar relations in 
Proposition~\ref{main-proposition}. 

The original statement of Theorem~\ref{intro:main-thm-surj} assumes that 
\(X\) contains a line defined over \(k\), which was used to prove 
Lemma~\ref{hyperplane-class}. The referee suggested a wearker assumption 
that \(F(X)\) contains a zero-cycle of degree \(1\), and pointed out 
that it is equivalent to the assumption that \(X\) contains a one-cycle 
of degree \(1\), see Lemma~\ref{lem:equiv-assump}.
  \begin{proposition}[=Proposition~\ref{main-proposition}] 
    \label{intro:main-prop}
    Let \(X\subset \mathbb{P}^{n+1}_k\) be a smooth cubic hypersurface 
    over a field \(k\), and let \(F(X)\) be the variety of lines. 
    Denote by \(h_X\in \mathrm{CH}^1(X)\) the hyperplane section class. 
    Assume \(X\) has a one-cycle of degree \(1\). 
    Then for a given algebraic cycle \(\Gamma\in \mathrm{CH}_r(X)\) of 
    dimension \(r>1\) with degree \(e \coloneqq h^r_X\cdot \Gamma\), there 
    exists two \((r-1)\)-cycles \(\gamma_1, \gamma_2\in \mathrm{CH}_{r-1}(F(X))\) 
    such that 
    \begin{align}
    &2\Gamma+q_*p^*\gamma_1 \in \mathbb{Z}\cdot h^{n-r}_X; \label{key-relation-1}\\
    &(2e-3)\Gamma+q_*p^*\gamma_2\in \mathbb{Z}\cdot h_X^{n-r}. \label{key-relation-2}
    \end{align}
\end{proposition}

These relations imply that any \(r\)-cycle on \(X\) lies in the image 
of the cylinder map module a multiple of the class \(h_X^{n-r}\). 
In Lemma~\ref{hyperplane-class}, we further show that \(h^{n-r}_X\) 
is also in the image of the cylinder map.
It is noteworthy that the formula~\eqref{key-relation-1} already appeared
in~\cite{Shen-uni-gen} by Shen, and Mboro has presented the 
formula~\eqref{key-relation-2} assuming \(\Gamma\) is a smooth closed 
subvariey in \(X\), see~\cite[Thm. 0.8]{Mboro:CH2}. The main improvement
of this proposition is to establish the second relation for arbitrary
\(r\)-cycles.

The proof of Proposition~\ref{intro:main-prop} is based on the 
birational geometry of the Hilbert scheme of two points on a cubic 
hypersurface, developed by Galkin-Shinder~\cite{GS14} and 
Voisin~\cite{uni-CH0}. We will review this in Section~\ref{sec:geom-Hilb-sq}.   

In Section~\ref{sec:IHC-ITC}, we apply Theorem~\ref{intro:main-thm-surj} 
to prove the integral Hodge conjecture and integral Tate conjecture 
for the variety of lines on a smooth cubic fourfold. The integral 
Hodge conjecture asks whether the integral Hodge classes on a given smooth 
complex projective variety are algebraic.
Unlike the well-known Hodge conjecture, examples~\cite{Grabowski:IHC-abelian-3fold,
Totaro:IHC-3folds-Kod0,Voisin:IHC-3folds} and counterexamples~\cite{AH-IHC,
Benoist-Ottem:IHC-3-folds,Kollar:Trento,Totaro:cycle-cobord} both exist 
for the integral Hodge conjecture, making it an interesting question 
to determine whether the conjecture holds for a special class of varieties.
 
The \(F(X)\) of a smooth complex cubic fourfold \(X\) is a 
hyperk\"ahler manifold of \(K3^{[2]}\)-type. Mongardi and Ottem
~\cite{Mongardi-Ottem:Hodge-Curve} recently showed that 
\((2n-2)\)-Hodge classes on any projective hyperk\"ahler \(n\)-fold 
of \(K3\)-type or of the generalized Kummer type are generated by 
rational curves, thus concludes the integral Hodge conjecture 
for curve classes on those hyperk\"ahler varieties
Using Theorem~\ref{intro:main-thm-surj} 
and the integral Hodge conjecture for cubic fourfolds by 
Voisin~\cite{Voisin:Hodge-conjecture}, we provide an alternative 
proof for the specific hyperk\"ahler fourfold \(F(X)\).
\begin{corollary}[=Theorem~\ref{int-Hodge-one-cycl}]
Let \(X\) be a smooth complex cubic fourfold. The integral Hodge 
conjecture holds for curve classes on the hyperk\"ahler variety \(F(X)\) 
\end{corollary}

The arithmetic analogue of the Hodge conjecture is the Tate conjecture. 
Suppose \(V\) is a smooh proper variety over a finitely generated
field \(k\). Let \(\bar{k}\) be the separable closure of \(k\) and 
let \(G_k\coloneqq \mathrm{Gal}(\bar{k}/k)\) be the absolute Galois group. 
Fix a prime \(\ell\) that is invertible in \(k\). The integral Tate 
conjecture for one-cycles on \(V\) states the following
\begin{conjecture}
The cycle class map 
\begin{equation}
  \label{lim-int-Tate}
\mathrm{cl}_{\bar{V}} \colon \mathrm{CH}_{1}(V_{\bar{k}})\underset{\mathbb{Z}}{\otimes} \mathbb{Z}_{\ell}
\to \underset{G_{k'}}{\varinjlim}~\mathrm{H}^{2\dim V-2}_{\et}(V_{\bar{k}}, \mathbb{Z}_{\ell}(\dim V-1))^{G_{k'}}.
\end{equation}
is surjective, where the direct limit is taken over all immediate finite 
field extensions \(k\subset k'\subset \bar{k}\).
\end{conjecture} 
This conjecture, formulated by Schoen~\cite{int-Tate}, differs from 
the ``naive'' integral Tate conjecture, where algebraic cycles are 
defined over the initial field \(k\). In particular, if \(k\) is 
a finite field, Schoen confirmed the conjecture under the assumption 
of the classical Tate conjecture for surfaces over any finite extension 
of \(k\), see~\cite[Thm. 0.5]{int-Tate}.

We verify Schoen's integral Tate conjecture for variety of lines
on a smooth cubic fourfold.
\begin{corollary}[=Theorem~\ref{ITC}]
Let \(k\) be a finitely generated field of characteristic different 
from $2$ and $3$. Let \(X\) be a smooth cubic fourfold over \(k\). 
Then the cycle class map~\eqref{lim-int-Tate} for one-cycles on \(F(X)\) 
is surjective.
\end{corollary}

\subsection*{Acknowledgements.} I owe my gradtitude to my advisor Mingmin Shen for his guidance on this work. 
We appreciate the anonymous referee for valuable comments and suggestion to improve our main statement.
We are grateful to Prof. Shimada for sharing his references. 
We would like to thank Prof. Colliot-Th\'el\`ene, Prof. Taelman and Xuanyu Pan for communications on early version of the manuscript. 
Thanks to Ren\'e Mboro for bringing my attention to his work. 

\section{The Hilbert square of a cubic hypersurface}
  \label{sec:geom-Hilb-sq}

Let \(k\) be a field. Let \(X\subset \mathbb{P}^{n+1}_k\) be a smooth 
cubic hypersurface of dimension \(n\). The variety \(F(X)\) of lines 
on \(X\) has been extensively studied through Altman-Kleiman's 
work~\cite{AK77}. We recall that \(F(X)\) is a smooth projective 
variety of pure dimension \(2n-4\). Moreover, it is geometrically 
connected if \(n\geq 3\). 
Consider the universal \(\mathbb{P}^1\)-bundle
\[
P\coloneqq \{([\ell], x)\in F(X)\times X~|~ x\in \ell\subset X\}
\]
with two projections \(p\colon P\to F(X)\) and \(q\colon P\to X\), and
the induced cylinder homomorphism
\[
P_*=q_*p^*\colon \mathrm{CH}_{r-1}(F(X))\to \mathrm{CH}_{r}(X), ~r\geq 1.
\]
We aim to prove the following two formulae.

\begin{proposition}
  \label{main-proposition}
  Let $h_X\in \mathrm{CH}^1(X)$ be the hyperplane section class. Suppose
  that \(X\) has a one-cycle of degree \(1\). Then for a given algebraic 
  cycle \(\Gamma\in \mathrm{CH}_r(X)\) of dimension $r>1$ and degree 
  \(e\coloneqq h^r_X\cdot \Gamma\), there exists two \((r-1)\)-cycles
  \(\gamma_1, \gamma_2\in \mathrm{CH}_{r-1}(F(X))\) satisfying 
    \begin{align}
    &2\Gamma+q_*p^*\gamma_1 \in \mathbb{Z}\cdot h^{n-r}_X; \label{eqs:relation-one}\\
    &(2e-3)\Gamma+q_*p^*\gamma_2\in \mathbb{Z}\cdot h_X^{n-r}. \label{eqs:relation-two}
    \end{align}
\end{proposition}
Let \(X^{[2]}\) be the Hilbert scheme of two points on \(X\), and let
\[
P_X=\{([\ell], x)\in \mathbb{G}(1, n+1)\times X |x\in \ell
\subset \mathbb{P}^{n+1}\}.
\]
be the incidence variety of lines in \(\mathbb{P}^{n+1}_k\) meeting 
with \(X\). Galkin and Shinder~\cite{GS14} constructed a birational map 
\[
    \Phi: X^{[2]}\dashrightarrow P_X
\]
as follows:
let \(\tau\in X^{[2]}\) represent two distinct points in \(X\) or a 
tangent vector supported on a closed point in \(X\). If the line 
\(\ell_{\tau}\) in \(\mathbb{P}^{n+1}\) generated by \(\tau\) is not 
contained in \(X\), one defines \(\Phi(\tau)=(\ell_{\tau},z)\in P_X\) 
where \(z\) is the unique residue point of the intersection 
\(\ell_{\tau}\cap X\).

Apparently \(\Phi\) is not defined on the set of points \(\tau\) such 
that \(\ell_{\tau}\subset X\). Let us denote \(P_2\) the indeterminacy 
of \(\Phi\). We can view \(P_2\) as the relative symmetric product of 
the \(\mathbb{P}^1\)-bundle \(p\colon P\to F(X)\). Then \(P_2\) is a
\(\mathbb{P}^2\)-bundle over \(F(X)\) because the fiber over any
\([\ell]\in F(X)\) is the symmetric product \(\ell^{(2)}\cong \mathbb{P}^2\). 
In~\cite{uni-CH0}, Voisin described the explicit resoultion of the birational 
map \(\Phi\) given by blowing up along \(P_2\).

\begin{proposition}\cite[Proposition 2.9]{uni-CH0}
    \label{resolution-map}
    \begin{enumerate}
    \item The birational map \(\Phi\) can be resolved by 
    the blowing up \(\tau\colon \widetilde{X^{[2]}}\to X^{[2]}\) 
    along the smooth center \(P_2\subset X^{[2]}\). 

    \item The induced map \(\widetilde{\Phi}\colon \widetilde{X^{[2]}}\to P_X\)
    identifies \(\widetilde{\Phi}\) with the blow up \(\widetilde{P_X}\)
    of \(P_X\) along the smooth center \(P\subset P_X\). 

    \item The exceptional divisors of the blow ups \(\widetilde{X^{[2]}}\)
    and \(\widetilde{P_X}\) are identified via the isomorphism
    \(\widetilde{X^{[2]}}\cong \widetilde{P_X}\).
    \end{enumerate}
\end{proposition}
\begin{proof}
The same statements in~\cite[Proposition 2.9]{uni-CH0} are over the 
complex numbers. But we note the same argument could apply to
an arbitrary field.
\end{proof}
For the later use, we collect the materials of the resolution in the 
following diagram
\begin{equation}
  \label{resolution-diagram}
    \begin{tikzcd}
    \mathcal{E} \ar[rr, "\pi_1"] \ar[dd, "\pi_2"'] \ar[rd, hook, "j"] && 
    P \ar[d, hook, "i_1"] \ar[r, "q"] & X \\    
    &\widetilde{X^{[2]}} \ar[r, "\widetilde{\Phi}"] \ar[d, "\tau"] & 
    P_X \ar[ru, "\pi_X"'] & \\
    P_2 \ar[r, hook, "i_2"] & X^{[2]} \ar[ru, dashrightarrow, "\Phi"'], &
    \end{tikzcd}
\end{equation}
where \(\mathcal{E}\) is the exceptional divisor. 

\begin{lemma}
  \label{lem:exp-inter}
Let \(\pi_{F*}: P_2\to F\) be the projection. For any algebraic cycle 
\(\Xi\in \mathrm{CH}_k(X^{[2]})\), the \((k-2)\)-cycle 
\(\gamma:=\pi_{F*}i_2^*\Xi\) in \(\mathrm{CH}_{k-2}(F(X))\) satisfies
    \begin{equation}
    \label{eqs:exp-inter}
    {\pi_X}_*\widetilde{\Phi}_*(\mathcal{E}\cdot \tau^*\Xi)=q_*p^*\gamma
    \end{equation} 
    in \(\mathrm{CH}_{k-1}(X)\).
\end{lemma}

\begin{proof}
    Note that \(\mathcal{E}\cdot \tau^*\Xi=j_*j^* \tau^*\Xi\). 
    The diagram~\eqref{resolution-diagram} implies that
    \[
    {\pi_X}_*\widetilde{\Phi}_*(\mathcal{E}\cdot \tau^*\Xi)
   ={\pi_X}_*\widetilde{\Phi}_*j_*j^* \tau^*\Xi
   ={\pi_X}_*{i_1}_*{\pi_1}_*{\pi_2}^*{i_2}^*\Xi
   =q_*{\pi_1}_*{\pi_2}^*{i_2}^*\Xi.
    \]
    The projection \(\pi_1 \colon \mathcal{E}\to P\) is a 
    \(\mathbb{P}^2\)-bundle since \(\mathrm{codim}(P, P_X)=3\). 
    Moreover, the exceptional divisor fits into the cartesian diagram
    \begin{equation}
      \label{flat-pul-pro-push}
      \begin{tikzcd}
      \mathcal{E} \ar[r, "\pi_1"] \ar[d, "\pi_2"'] & P \ar[d, "p"]\\
      P_2 \ar[r, "\pi_F"] &F(X).
      \end{tikzcd}
    \end{equation}
    Hence \({\pi_1}_*{\pi_2}^*{i_2}^*\Xi=p^*{\pi_F}_*{i_2}^*\Xi\).
    Denote by \(\gamma\) the cycle \({\pi_F}_*{i_2}^*\Xi\in \mathrm{CH}_{k-2}(F(X))\). 
    It follows that
      \[
      {\pi_X}_*\widetilde{\Phi}_*(\mathcal{E}\cdot \tau^*\Xi)=q_*p^*\gamma.
      \]
\end{proof}

It is crucial to describe the divisor class of 
\(\mathcal{E}\in \mathrm{Pic}(\widetilde{X^{[2]}})\). For this purpose,
we shall consider the natural morphism
\[
\varphi\colon X^{[2]}\to \mathbb{G}(1, n+1),
\]
which assigns to \(\tau\in X^{[2]}\) the generating line 
\(\ell_{\tau}\subset \mathbb{P}^{n+1}\). Denote by \(\mathcal{U}\) 
(resp. \(\mathbb{P}(\mathcal{U})\)) the tautological rank \(2\) 
subbundle (resp. \(\mathbb{P}^1\)-bundle) over \(\mathbb{G}(1, n+1)\).
The pullback \(\mathcal{Q}:=\mathbb{P}(\varphi^*\mathcal{U})\) 
is a \(\mathbb{P}^1\)-bundle over \(X^{[2]}\). Consider the composition 
\[
\alpha\colon \mathcal{Q}\to \mathbb{P}(\mathcal{U})\overset{\pi_P}{\longrightarrow} \mathbb{P}^{n+1}  
\]
where \(\pi_P\) is the natural projection. In the proof of Proposition
~\ref{resolution-map}, Voisin showed that the divisor \(\alpha^{-1}(X)\) 
in \(\mathcal{Q}\) consists of two components \(\widetilde{X^{[2]}}\) 
and the blow up \(\widetilde{X\times X}\) of \(X\times X\) along the 
diagonal \(\Delta_X\).
Let 
\[
\rho\colon \widetilde{X\times X}\to X\times X
\] 
denote the blowing up, and \(E_{\Delta,X}\) denote the exceptional
divisor. The involution \((x,y)\mapsto (y,x)\) on \(X\times X\) induces
a natural involution on \(\widetilde{X\times X}\), and the quotient by 
the involution yields to a double covering  
\[
\sigma\colon \widetilde{X\times X}\to X^{[2]}.
\]
We have a divisor class \(\delta\in \mathrm{Pic}(X^{[2]})\) such that 
\(\sigma^*\delta=E_{\Delta, X}\) in \(\mathrm{Pic}(\widetilde{X\times X})\).
The explicit definition is \(\delta:=-c_1(\sigma_*\mathcal{O}_{\widetilde{X\times X}}/\mathcal{O}_{X^{[2]}})\), see~\cite[Eq. (12)]{Shen-uni-gen}. We call 
\(\delta\) the \emph{half diagonal class} because \(2\delta\) is 
equivalent to the diagonal class \(E_{\Delta, X}\) on \(X^{[2]}\).
Given algebraic cycles \(\alpha, \beta\in \mathrm{CH}^*(X)\), we write
\[
\alpha\hat{\otimes} \beta \coloneqq \sigma_*\rho^*(\alpha\otimes \beta)
\in \mathrm{CH}^*(X^{[2]}). 
\] 

\begin{lemma}
  \label{expdiv}
  Denote by \(h\in \mathrm{CH}^1(\mathbb{P}^{n+1})\) the hyperplane 
  class, and by \(h_X\in \mathrm{CH}^1(X)\) the restriction of \(h\) 
  to \(X\). Let \(h_{\mathcal{Q}}\in \mathrm{Pic}(\mathcal{Q})\) be 
  the divisor class \(\alpha^*h_X\). 
  Then we have
\[
\mathcal{E}=-h_{\mathcal{Q}}|_{\widetilde{X^{[2]}}}+\tau^*(2h_X\hat{\otimes}1-3\delta)
\in \mathrm{Pic}(\widetilde{X^{[2]}}),
\]
and
\[
c_1(\varphi^*\mathcal{U})=-h_X\hat{\otimes}1+\delta\in \mathrm{Pic}(X^{[2]}).
\]
\end{lemma}
\begin{proof}
For the first equation, the sign before the second term is negative in~\cite[Lemma 4.3]{Shen-uni-gen}.
But the proof in loc. cit. shows that it should be positive. The second equation
can be found at the end of the proof of~\cite[Lemma 4.3]{Shen-uni-gen}.
\end{proof}

To obtain the relations~\eqref{eqs:relation-one} and~\eqref{eqs:relation-two}
for a given cycle \(\Gamma\) on \(X\), we shall construct certain cycles
\(\Xi\) on \(X^{[2]}\) and reinterpret \(\Gamma\) via the equation
~\eqref{eqs:exp-inter}, which involves some technical computations. 
Let us set the following cartesian diagram to track notations in the 
coming proofs
\begin{equation}
  \label{eqs:cartesqs}
    \begin{tikzcd}
    \widetilde{X\times X}\coprod \widetilde{X^{[2]}} \ar[r, "\Psi\coprod \widetilde{\Phi}"] 
    \ar[d] & P_X \ar[d, hook, "i'"] \ar[r, "\pi_X"] & X \ar[d, hook, "i_X"]\\
    \mathcal{Q}\ar[r] \ar[d] & 
    \mathbb{P}(\mathcal{U}) \ar[d, "\pi_G"] \ar[r, "\pi_P"]& \mathbb{P}^{n+1}\\
    X^{[2]}\ar[r, "\varphi"] &\mathbb{G}(1, n+1).
    \end{tikzcd}
\end{equation}
The map \(\Psi\colon \widetilde{X\times X}\to P_X\) given by
\(\Psi(x,y)=(x, \ell_{xy})\) fits into the commutative diagram
  \begin{equation} 
    \label{dobcov}
    \begin{tikzcd}
    X^{[2]} & \widetilde{X\times X} \ar[r, "\Psi"] \ar[d, "\rho"] \ar[l, "\sigma"'] & P_X \ar[d,"\pi_X"]\\
    &X\times X\ar[r, "p_1"] &X,
    \end{tikzcd}
  \end{equation}
where \(p_1\) is the projection to the first factor.

\begin{lemma}
  \label{lem:cycle-correspondence}
  With the notations in diagrams~\eqref{resolution-diagram},
  \eqref{eqs:cartesqs} and~\eqref{dobcov}, given any algebraic 
  cycle \(Z\in \mathrm{CH}^i(X^{[2]})\), we have
    \begin{equation}
    \label{eqs:comm-correspond}
    {\pi_X}_*\widetilde{\Phi}_*\tau^*Z+{\pi_X}_*\Psi_*\sigma^*Z=i_X^*{\pi_P}_*{\pi_G}^*\varphi_*Z.
    \end{equation}
\end{lemma}

\begin{proof}    
    Regard \(\widetilde{X^{[2]}}\) and \(\widetilde{X\times X}\) 
    as cycle correspondences on \(X^{[2]}\times P_X\). As seen from
    diagrams~\eqref{resolution-diagram} and~\eqref{dobcov}, the operator 
    \(\widetilde{\Phi}_*\circ  \tau^*+\Psi_*\circ \sigma^*\) on the
    Chow groups is inducded by the correspondence 
    \(\widetilde{X^{[2]}}+\widetilde{X\times X}\).
    It follows from the cartesian diagram~\eqref{eqs:cartesqs} that
    \[
    \widetilde{\Phi}_*\circ  \tau^*+\Psi_*\circ \sigma^* = {i'}^*\circ \pi^*_G\circ \varphi_*.
    \]  
    The assertion of the lemma follows by composing the pushforward
    \({\pi_X}_*\) on both sides and \(i_X^*{\pi_P}_* = {\pi_X}_*{i'}^*\).
\end{proof}

Suppose that \(\eta\) is a one-cycle on \(X\), and \(\Gamma\subset X\)
is a closed subvariety of dimension \(r\). Then 
\([\Gamma]\otimes \eta + \eta\otimes [\Gamma]\) and \([\Gamma\times \Gamma]\) 
are invariant cycles on \(X\times X\) under the involution \((x,y)\mapsto (y,x)\). 
By~\cite[Corollary 2.4]{uni-CH0} there exist two cycles \(\Sigma_1\) and 
\(\Sigma_2\) of \(X^{[2]}\) satisfying
\begin{equation}
   \label{eqs:symcyc}
\rho_*\sigma^*\Sigma_1 = [\Gamma]\otimes \eta + \eta\otimes [\Gamma],~ 
\rho_*\sigma^*\Sigma_2 = [\Gamma\times \Gamma].
\end{equation}
More precisely, we set 
\begin{align}
\Sigma_1\colon &= \sigma_*\rho^*([\Gamma]\otimes \eta)\in \mathrm{CH}_{r+1}(X^{[2]}); \label{eqs:sym-cycl-1}\\
\Sigma_2\colon &= [\sigma(\widetilde{\Gamma\times \Gamma})]\in \mathrm{CH}_{2r}(X^{[2]}), \label{eqs:sym-cycl-2}
\end{align}
where \(\widetilde{\Gamma\times \Gamma}\) is the strict transform of 
\(\Gamma\times \Gamma\) in the blow up \(\widetilde{X\times X}\).
We remark that \([\sigma(\widetilde{\Gamma\times \Gamma})]\) 
indicates the cycle of the closed image with reduced structure rather than the pushforward 
\(\sigma_*[\widetilde{\Gamma\times \Gamma}]\). The latter is twice of 
the former.

Let us check the relations~\eqref{eqs:symcyc}. It is direct to see
\[
\rho_*\sigma^*\Sigma_1=\rho_*\sigma^*\sigma_*\rho^*([\Gamma]\otimes \eta)
=\rho_*\rho^*([\Gamma]\otimes \eta+\eta\otimes [\Gamma])
=[\Gamma]\otimes \eta + \eta\otimes [\Gamma].
\]
Due to the flatness of \(\sigma\) we have \(\sigma^*[\sigma(\widetilde{\Gamma\times \Gamma})]
=[\sigma^{-1}(\sigma(\widetilde{\Gamma\times \Gamma}))]
=[\widetilde{\Gamma\times \Gamma}]\). Hence
\[
\rho_*\sigma^*\Sigma_2=
\rho_*[\widetilde{\Gamma\times \Gamma}]=[\Gamma\times\Gamma].
\]

The next lemma includes necessary computations for the proof
of the Proposition~\ref{main-proposition}.

\begin{lemma}
  \label{cycle-symm-prod}
  Let \(X\) be a smooth cubic hypersurface, \(h_X\) be the hyperplane
  section class, and \(\delta\) be the half diagonal class on \(X^{[2]}\). 
  Let \(\Gamma\subset X\) be a closed subvariety of dimension \(r>1\)
  with degree \(e:=\Gamma\cdot h_X^{n-r}\). Let \(\eta\in \mathrm{CH}_1(X)\) 
  be a one-cycle of degree \(1\). Consider the algebraic 
  cycles \(\Sigma_1\) and \(\Sigma_2\) constructed  
  in~\eqref{eqs:sym-cycl-1} and~\eqref{eqs:sym-cycl-2}. 
\begin{enumerate}
\item With the notations in diagram~\eqref{dobcov}, we show
\begin{align*}
{\pi_X}_*\Psi_*\sigma^*\Sigma_1&=0,\\
{\pi_X}_*\Psi_*\sigma^*(h_X\hat{\otimes}1\cdot \Sigma_1)&=\Gamma,\\
{\pi_X}_*\Psi_*\sigma^*(\delta\cdot \Sigma_1)&=0.
\end{align*}
\item With the notations in diagram~\eqref{dobcov}, we show 
\[
{\pi_X}_*\Psi_*\sigma^*(\Sigma_2\cdot (h_X\hat{\otimes} 1)^k\cdot \delta^{l})=
\begin{cases}
0, & \forall ~ 0\leq k, l \leq r-1,\\
e\cdot \Gamma, & k=r, l=0,\\
(-1)^{r+1}\Gamma, & k=0, l=r.
\end{cases}
\]
\end{enumerate} 
\end{lemma}

\begin{proof}
By the diagram~\eqref{dobcov} we have \({\pi_X}_*\Psi_*={p_1}_*\rho_*\).
It follows from~\eqref{eqs:symcyc} that 
\[
{\pi_X}_*\Psi_*\sigma^*\Sigma_1={p_1}_*\rho_*\sigma^*\Sigma_1
={p_1}_*(\Gamma\otimes \eta+\eta\otimes \Gamma)=0.
\]

Recall that \(\sigma^*(\alpha\hat{\otimes} \beta) = \sigma^*\sigma_*\rho^*(\alpha\otimes \beta)
=\rho^*(\alpha\otimes \beta+\beta\otimes \alpha)\).
Therefore 
  \begin{align*}
    {\pi_X}_*\Psi_*\sigma^*(h_X\hat{\otimes}1\cdot\Sigma_1)
    &={p_1}_*\rho_*(\sigma^*(h_X\hat{\otimes}1)\cdot \sigma^*\Sigma_1)\\
    &={p_1}_*\rho_*(\rho^*(h_X\otimes 1+1\otimes h_X)\cdot \sigma^*\Sigma_1)\\
(\textrm{projection formula})&={p_1}_*((h_X\otimes 1+1\otimes h_X)\cdot \rho_*\sigma^*\Sigma_1)\\
    &={p_1}_*((h_X\otimes 1+1\otimes h_X)\cdot (\Gamma\otimes \eta+\eta\otimes \Gamma)).
  \end{align*}
Note that \(\dim (h_X\cdot \Gamma) > 0\) since \(\dim \Gamma=r>1\). 
Hence the last pushforward equals to \(\deg(h_X\cdot \eta)\cdot \Gamma=\Gamma\).

To prove \({\pi_X}_*\Psi_*\sigma^*(\delta\cdot \Sigma_1)=0\) we shall 
use the diagram
\begin{equation}
  \label{diag-blow-up}
  \begin{tikzcd}
  E_{\Delta, X} \ar[r, hook, "j_E"] \ar[d, "\pi_\Delta"] 
  &\widetilde{X\times X} \ar[r, "\sigma"] \ar[d, "\rho"]
  &X^{[2]}\\
  X \ar[r, hook, "\iota_{\Delta}"] & X\times X \ar[ur, dashrightarrow, "\mu"]
  \end{tikzcd} 
\end{equation}
where \(\iota_{\Delta}\) is the diagonal embedding. Then we have
    \[
      {\pi_X}_*\Psi_*\sigma^*(\delta \cdot \Sigma_1)
      ={\pi_X}_*\Psi_*(E_{\Delta,X}\cdot \sigma^*\Sigma_1)
      ={p_1}_*\rho_*({j_E}_*{j_E}^*\sigma^*\Sigma_1)
    \]
It follows from \(\rho\circ {j_E} = \iota_\Delta\circ \pi_\Delta\) and 
\(p_1\circ \iota_\Delta=\operatorname{id}_X\) that
    \begin{align*}
      {p_1}_*\rho_*({j_E}_*{j_E}^*\sigma^*\Sigma_1)
     &={\pi_\Delta}_*{j_E}^*\sigma^*\Sigma_1\\
     &={\pi_\Delta}_*{j_E}^*\rho^*(\Gamma\otimes \eta+\eta\otimes \Gamma)\\
     &={\pi_\Delta}_*{\pi_\Delta}^*{\iota_\Delta}^*(\Gamma\otimes \eta+\eta\otimes \Gamma)\\
     &=2{\pi_\Delta}_*{\pi_\Delta}^*(\Gamma\cdot \eta).
    \end{align*}
Note that \(\pi_\Delta\colon E_{\Delta,X}\to X\) is a projective bundle of
positive relative dimension.
Hence \({\pi_\Delta}_*{\pi_\Delta}^*=0\). Therefore assertion \((1)\) 
is complete.

For the assertion \((2)\), we first note 
\[
\sigma^*(h_X\hat{\otimes} 1)^k=\rho^*(h_X\otimes 1+1\otimes h_X)^k
=\sum_{s=0}^{k}\binom{k}{s}\rho^*(h_X^s\otimes h_X^{k-s}).
\]
Recall that the cycle class \(\sigma^*\Sigma_2\) is represented by 
\([\widetilde{\Gamma\times\Gamma}]\). Similar argument as above
yields
\[
  {\pi_X}_*\Psi_*\sigma^*(\Sigma_2\cdot (h_X\hat{\otimes} 1)^k\cdot \delta^{l})
  ={p_1}_*(\rho_*([\widetilde{\Gamma\times\Gamma}]\cdot E^l_{\Delta,X})\cdot \sum_{s=0}^{k}\binom{k}{s}h_X^s\otimes h_X^{k-s}).
\]
Let us carry out computations case by case.
\begin{itemize}
\item Suppose that \(l=0, k\leq r\). It follows from 
\(\rho_*\sigma^*\Sigma_2=[\Gamma\times \Gamma]\) that
\begin{align*}
{\pi_X}_*\Psi_*\sigma^*(\Sigma_2\cdot (h_X\hat{\otimes}1)^k)
=&{p_1}_*([\Gamma\times \Gamma]\cdot \sum_{s=0}^{k}\binom{k}{s}h_X^s\otimes h_X^{k-s})\\
=&\sum_{s=0}^{k}{{k}\choose{s}}{p_1}_*(h_X^s\cdot \Gamma\otimes h_X^{k-s}\cdot \Gamma)\\
=&\begin{cases}
0, & k< r;\\
\deg \Gamma \cdot \Gamma, & k=r.
\end{cases}
\end{align*}
We see \({p_1}_*(h_X^s\cdot \Gamma\otimes h_X^{k-s}\cdot \Gamma)\)
is nontrivial if and only if \(k=r\) and \(s=0\) since the dimension of 
\(\Gamma\) is \(r\).
\item Suppose that $1\leq l\leq r-1$. We claim that
\[
  \rho_*([\widetilde{\Gamma\times\Gamma}]\cdot E^l_{\Delta,X})=0.
\]

Consider the excess normal bundle \(Q\) of the exceptional divisor
\(E_{\Delta,X}\). It is the quotient bundle given by the exact sequence
\begin{equation}
  \label{excess-normal-bdl}
0\rightarrow \mathcal{O}_E(-1)\rightarrow \pi^*_{\Delta}N_{X/X\times X} \rightarrow Q\rightarrow 0
\end{equation}
where \(\mathcal{O}_E(1)\) is tautological line bundle. 
Fulton's blow-up formula~\cite[Thm. 6.7]{Fulton} associated with the 
blow-up diagram~\eqref{diag-blow-up} presents the cycle class of
the strict transform \(\widetilde{\Gamma\times\Gamma}\) as follows
\[
[\widetilde{\Gamma\times\Gamma}]=\rho^*(\Gamma\times \Gamma)-{j_E}_*\{c(Q)\cap \pi^*_{\Delta}s(\Gamma, \Gamma\times \Gamma)\}_{2r}
\]
where $s(\Gamma, \Gamma\times \Gamma)$ is the Segre class of the 
closed subvariety $\Gamma\subset \Gamma\times \Gamma$, and the notation
\(\{\alpha\}_{2r}\) indicates the \(2r\)-th dimensional part of 
a given total cycle class \(\alpha\in \bigoplus_{k\geq 0} \mathrm{CH}_k\).
Using the blow-up formula, the projection formulae for the morphism
\(\rho\) and \(j_E\), and \({\pi_{\Delta}}\circ {\iota_\Delta}=\rho\circ j_E\), 
the class \(\rho_*([\widetilde{\Gamma\times\Gamma}]\cdot E_{\Delta}^{l})\)
equals to
\begin{equation}
\label{eqs:strict-transform-inter}
[\Gamma\times \Gamma]\cdot \rho_*E_{\Delta, X}^{l}
-{\iota_\Delta}_*{\pi_{\Delta}}_*(\{c(Q)\cdot \pi^*_{\Delta}s(\Gamma, \Gamma\times \Gamma)\}_{2r}\cdot j^*_E E^l_{\Delta, X}).
\end{equation}
For the first term of~\eqref{eqs:strict-transform-inter} we note that
\[
\rho_*E_{\Delta, X}^l=\rho_*{j_E}_*{j^*_E}E_{\Delta, X}^{l-1}
={\iota_\Delta}_*{\pi_\Delta}_*c_1(\mathcal{O}_E(-1))^{l-1}.
\]
The class \({\pi_\Delta}_*c_1(\mathcal{O}_E(-1))^{l-1}\)
is zero, because \(l\) is less than the rank of the normal bundle 
\(N_{X/X\times X}\). Hence the first term of~\eqref{eqs:strict-transform-inter} 
vanishes.

It remains to show the second term of~\eqref{eqs:strict-transform-inter} 
is zero. We simply denote by \(s_i\) the \(i\)-th Segre class of 
\(s(\Gamma, \Gamma\times \Gamma)\), which is a \((r-i)\)-cycle supported 
on \(\Gamma\). Then the \(2r\)-th dimensional part of the total class 
\(c(Q)\cdot \pi^*_{\Delta}s(\Gamma, \Gamma\times \Gamma)\) is the sum
\[
\sum_{i+t=n-r-1} c_t(Q)\cdot {\pi_\Delta}^*s_i\in \mathrm{CH}_{2r}(E_{\Delta, X}).
\]
By the projection formula
\({\pi_{\Delta}}_*(c_t(Q)\cdot {\pi_\Delta}^*s_i\cdot j^*_E E^l_{\Delta, X})\)
is equal to
\[
s_i\cdot {\pi_\Delta}_*(c_t(Q)\cdot c_1(\mathcal{O}_E(-1))^l).
\]
Denote by \(\xi\) the first Chern class \(c_1(\mathcal{O}_E(1))\), 
and by \(N\) the normal bundle \(N_{X/X\times X}\). Using the exact 
sequence~\eqref{excess-normal-bdl}, the total Chern class of the 
excess normal bundle \(Q\) equals to
\[
c(Q)=\frac{\pi^*_{\Delta}c(N)}{c(\mathcal{O}_E(-1))}=\frac{\pi^*_{\Delta}c(N)}{1-\xi}
=\pi^*_{\Delta}c(N)\cdot \sum_{i\geq 0}^{\infty} \xi^i,
\] 
which implies that the \(t\)-th term 
\(c_t(Q)=\displaystyle{\sum^{t}_{j=0}} \pi^*_{\Delta}c_j(N)\cdot \xi^{t-j}\). 
As a consequence, we have
\begin{equation}
   \label{Chernclass-quot-bdl}
{\pi_\Delta}_*(c_t(Q)c_1(\mathcal{O}_E(-1))^l))=(-1)^l\sum_{j=0}^{t} c_j(N)\cdot {\pi_{\Delta}}_*\xi^{l+t-j}. 
\end{equation}
Note that \(t\leq n-r-1\) and \(l\leq r-1\). Hence \(l+t-j < n-1\) for all 
\(0\leq j\leq t\), which asserts \({\pi_\Delta}_*\xi^{l+t-j}=0\). 
Therefore the second term of~\eqref{eqs:strict-transform-inter} is also zero. 
Hence \(\rho_*([\widetilde{\Gamma\times\Gamma}]\cdot E_{\Delta}^{l})=0\)

\item Suppose that $k=0, l=r$ and $r<n$. We need to compute the 
cycle~\eqref{eqs:strict-transform-inter} for \(l=r\). 
By the same argument as bove, we have \(\rho_*{E^r_{\Delta,X}}=0\)
since \(r < n= \mathrm{rk}(N_{X/X\times X})\). By the previous computation, 
the second term of~\eqref{eqs:strict-transform-inter} for \(l=r\) is
\[
-\sum_{i+t=n-r-1}{\iota_\Delta}_*(s_i\cdot {\pi_\Delta}_*(c_t(Q)c_1(\mathcal{O}_E(-1))^r)).
\]
We have shown via~\eqref{Chernclass-quot-bdl} that 
\({\pi_\Delta}_*\xi^{r+t-j}\) is nonzero only if \(j=0, t=n-r-1\), 
and \(i=0\). Since the \(0\)-th Segre class \(s_0=[\Gamma]\),
hence
\begin{align*}
{\pi_X}_*\Psi_*\sigma^*(\Sigma_2\cdot \delta^r)
&=(-1)^{r+1}{p_1}_*{\iota_\Delta}_*(s_0\cdot c_0(N)\cdot {\pi_\Delta}_*c_1(\mathcal{O}_E(1))^{n-1}))\\
&=(-1)^{r+1}s_0\cap [X]\\
&=(-1)^{r+1}[\Gamma].
\end{align*}

When $k=0, l=r=n$, the variety \(\Gamma\) is the total space \(X\). 
Then the cycle \(\Sigma_2\) is total class \([X^{[2]}]\), and 
\(\sigma^*\Sigma_2=[\widetilde{X\times X}]\).
It is straightforward to see that
\begin{align*}
{\pi_X}_*\Psi_*\sigma^*(\Sigma_2\cdot \delta^r)
&={p_1}_*\rho_*([\widetilde{X\times X}]\cdot E_{\Delta, X}^{n})\\
&={\pi_\Delta}_*{j_E}_*{j^*_E}E^{n-1}_{\Delta, X}\\
&={\pi_\Delta}_*c_1(\mathcal{O}_E(-1))^{n-1}\\
&=(-1)^{n+1}[X].
\end{align*}
\end{itemize}
\end{proof}

\begin{remark}\label{rk1}
Over a non-closed field, a cubic hypersurface does not necessarily contain
one-cycles of degree \(1\), see \cite{lines-cubic-finite}.
\end{remark}

\section{Proof of surjectivity of the cylinder maps}

\begin{proof}[Proof of Proposition~\ref{main-proposition}]
Suppose that \(\Gamma\subset X\) is an irreducible closed subvariety
of dimension \(r > 1\) and degree \(e=\Gamma\cdot h_X^{n-r}\).
Let \(\Sigma_1\) be the algebraic cycle given in~\eqref{eqs:sym-cycl-1}. 
By Lemma~\ref{lem:exp-inter} 
\[
{\pi_X}_*\widetilde{\Phi}_*(\mathcal{E}\cdot \tau^*\Sigma_1)=q_*p^*\gamma_1,
\] 
where \(\gamma_1={\pi_F}_*{i_2}^*\Sigma_1\in \mathrm{CH}_{r-1}(F(X))\).
Recall from Lemma~\ref{expdiv} that the exceptional divisor
\(\mathcal{E}=-h_{\mathcal{Q}}|_{\widetilde{X^{[2]}}}+
\tau^*(2h_X\hat{\otimes}1-3\delta).\)
Hence we have
\[ 
{\pi_X}_*\widetilde{\Phi}_*(\mathcal{E}\cdot \tau^*\Sigma_1)
=-{\pi_X}_*\widetilde{\Phi}_*(h_{\mathcal{Q}}|_{\widetilde{X^{[2]}}}\cdot\tau^*\Sigma_1)
+{\pi_X}_*\widetilde{\Phi}_*(\tau^*(2h_X\hat{\otimes}1-3\delta)\cdot\tau^*\Sigma_1).
\]
Note that \(h_{\mathcal{Q}}=\alpha^*h\). We see from the commutative 
diagram~\eqref{eqs:cartesqs} that  
\[
{\pi_X}_*\widetilde{\Phi}_*(h_{\mathcal{Q}}|_{\widetilde{X^{[2]}}}\cdot\tau^*\Sigma_1)
=h_X\cdot {\pi_X}_*\widetilde{\Phi}_*\tau^*\Sigma_1.
\]
It follows from Lemma~\ref{lem:cycle-correspondence} and Lemma
~\ref{cycle-symm-prod} that
\begin{align*}
{\pi_X}_*\widetilde{\Phi}_*\tau^*\Sigma_1
&=i_X^*{\pi_P}_*{\pi_G}^*\varphi_*\Sigma_1-\pi_*\Psi_*\sigma^*\Sigma_1\\
&=i_X^*{\pi_P}_*{\pi_G}^*\varphi_*\Sigma_1.
\end{align*}
Note that 
\({\pi_P}_*{\pi_G}^*\varphi_*\Sigma_1\in \mathrm{CH}^{n-r-1}(\mathbb{P}^{n+1})\).
Hence \({\pi_X}_*\widetilde{\Phi}_*(h_{\mathcal{Q}}|_{\widetilde{X^{[2]}}}\cdot\tau^*\Sigma_1)\) 
must be a multiple of the class \(h_X^{n-r}\). With the same process 
we have
\begin{align*}
&{\pi_X}_*\widetilde{\Phi}_*(\tau^*(2h_X\hat{\otimes}1-3\delta)\cdot\tau^*\Sigma_1)
\\
=&i_X^*{\pi_P}_*{\pi_G}^*\varphi_*((2h_X\hat{\otimes}1-3\delta)\cdot\Sigma_1)
-{\pi_X}_*\Psi_*\sigma^*((2h_X\hat{\otimes}1-3\delta)\cdot\Sigma_1)\\
\equiv& -\pi_*\Psi_*\sigma^*((2h_X\hat{\otimes}1-3\delta)\cdot\Sigma_1)~\mathrm{mod}~ \mathbb{Z}\cdot h_X^{n-r}
\end{align*}
Again Lemma~\ref{cycle-symm-prod} yields
\(\pi_*\Psi_*\sigma^*((2h_X\hat{\otimes}1-3\delta)\cdot \Sigma_1)=2\Gamma\). 
Therefore the class \(2\Gamma+q_*p^*\gamma_1\) is a multiple of \(h^{n-r}_X\).

Let \(\Sigma_2\) be the second algebraic cycle in~\eqref{eqs:sym-cycl-2}. 
To simplify the notations, let \(g\) denote the divisor class 
\(c_1(\varphi^*\mathcal{U})=-h_X\hat{\otimes}1+\delta\) in Lemma~\ref{expdiv}. 
The second formula~\eqref{key-relation-2} will be derived by comupting 
the following cycle
\[
{\pi_X}_*\widetilde{\Phi}_*(\mathcal{E}\cdot \tau^*(\Sigma_2\cdot g^{r-1})).
\] 
We see \(\dim (\Sigma_2\cdot g^{r-1})=r+1\). Again by 
Lemma~\ref{lem:exp-inter}
\begin{equation}
   \label{second-key-cycle}
\pi_*\widetilde{\Phi}_*(\mathcal{E}\cdot \tau^*(\Sigma_2\cdot g^{r-1}))=q_*p^*\gamma_2,
\end{equation}
where \(\gamma_2={\pi_F}_*{i_2}^*(\Sigma_2\cdot g^{r-1})\in \mathrm{CH}_{r-1}(F(X))\).
The same argument as above yields
\[
{\pi_X}_*\widetilde{\Phi}_*(\mathcal{E}\cdot \tau^*(\Sigma_2\cdot g^{r-1}))
=-h_X\cdot {\pi_X}_*\widetilde{\Phi}_*\tau^*(\Sigma_2\cdot g^{r-1})
+{\pi_X}_*\widetilde{\Phi}_*\tau^*((2h_X\hat{\otimes}1-3\delta)\cdot \Sigma_2\cdot g^{r-1}).
\]
It follows from Lemma~\ref{lem:cycle-correspondence} that
\[
{\pi_X}_*\widetilde{\Phi}_*\tau^*(\Sigma_2\cdot g^{r-1})=i_X^*{\pi_P}_*{\pi_G}^*\varphi_*(\Sigma_2\cdot g^{r-1})-{\pi_X}_*\Psi_*\sigma^*(\Sigma_2\cdot g^{r-1}).
\]
The first term on the right side is a multiple of the class 
\(h_X^{n-r-1}\). The second term is equal to
\[
{\pi_X}_*\Psi_*\sigma^*(\Sigma_2\cdot \sum_{k=0}^{r-1}\binom{r-1}{k}(-h_X\hat{\otimes}1)^k\cdot\delta^{r-k-1}).
\]
This class vanishes as a result of the assertion \((2)\) of 
Lemma~\ref{cycle-symm-prod}. Therefore the class 
\(h_X\cdot {\pi_X}_*\widetilde{\Phi}_*\tau^*(\Sigma_2\cdot g^{r-1})\)
is a multiple of \(h_X^{n-r}\). By the same argument we can see
\({\pi_X}_*\widetilde{\Phi}_*\tau^*((2h_X\hat{\otimes}1-3\delta)\cdot \Sigma_2\cdot g^{r-1})\) equals to
\[
-{\pi_X}_*\Psi_*\sigma^*((2h_X\hat{\otimes}1-3\delta)\cdot \Sigma_2\cdot \sum_{k=0}^{r-1}
\binom{r-1}{k}(-h_X\hat{\otimes}1)^k\cdot\delta^{r-k-1})
\]
modulo a multiple of \(h_X^{n-r}\). The assertion \((2)\) of Lemma
~\ref{cycle-symm-prod} shows that the above class is 
\[
(-1)^{r}2\deg \Gamma\cdot \Gamma+(-1)^{r+1}3\Gamma.
\]
To conclude the formula~\eqref{eqs:relation-two} we may replace 
\(\gamma_2\) by \((-1)^{r-1}\gamma_2\). Then the above computation
implies that 
\[
(2\deg\Gamma-3)\Gamma+ q_*p^*\gamma_2 \in \mathbb{Z}\cdot h_X^{n-r}.
\]
Through the linear combination, one easily extend relations
~\eqref{eqs:relation-one} and~\eqref{eqs:relation-two} to
cycles.
\end{proof}
The following corollary is immediate.
\begin{corollary}
  \label{surj-mod}
  Let \(X\subset \mathbb{P}^{n+1}_k\) be a smooth cubic hypersurface
  of dimension $n\geq 3$ over a field $k$, and let $F(X)$ be the 
  variety of lines on \(X\). Assume that \(X\) has a one-cycle of
  degree one. Fix any \(r\geq 1\), the cylinder map
    \[
    q_*p^*: \mathrm{CH}_{r-1}(F(X))\rightarrow \mathrm{CH}_r(X)
    \]
    is surjective modulo the subgroup generated by the class
    \(h^{n-r}_X\). 
\end{corollary}

\begin{proof} 
The conclusion for \(r=1\) had been proved by 
Shen~\cite[Proposition 4.2]{Shen-uni-gen}.
We will prove the case for \(r > 1\). Suppose that 
\(\Gamma=\sum_{i} n_i\Gamma_i\) is a \(r\)-cycle on \(X\) with 
irreducible components \(\Gamma_i\). Let \(e_i\) be the degree of each 
\(\Gamma_i\). Let \(\gamma_{i1}\) and \(\gamma_{i2}\) be the 
\((r-1)\)-cycles associated with \(\Gamma_i\) satisfying the 
two formulae in Proposition~\ref{main-proposition}. It follows that
\[
\Gamma_i+q_*p^*((e_i-1)\gamma_{i1}-\gamma_{i2})\in \mathbb{Z}\cdot h_X^{n-r}.
\]
By linear combination \(\Gamma\) is contained in the subgroup generated by 
\(q_*p^*\mathrm{CH}_{r-1}(F(X))\) and \(\mathbb{Z}\cdot h^{n-r}_X\).
\end{proof}
It remains to show that the class of hyperplane intersections 
\(h^{n-r}_X\) is also contained in the image of the cylinder map. 
Following the idea in the proof of Corollary~\ref{surj-mod}, we shall seek for 
a pair of coprime integers \((a,b)\) and corresponding
\((r-1)\)-cycles \(\alpha_{r-1}\) and \(\beta_{r-1}\) satisfying 
\begin{equation}
\label{coprime-hyperplane-class}
q_*p^*\alpha_{r-1}=a h^{n-r}_X,~q_*p^*\beta_{r-1}=b h^{n-r}_X,
\end{equation}
for which the next Lemma aims to prove.

\begin{lemma}
  \label{hyperplane-class}
  Let \(X\) be a smooth hypersurface of dimension \(n\geq 3\) over
  a field \(k\). Let \(h_X\) be the hyperplane section class on \(X\). 
  Assume that \(F(X)\) contains a zero-cycle of degree \(1\). Then 
  for any \(1\leq r < n-1\), there exists cycles \(\alpha_{r-1}\) 
  and \(\beta_{r-1}\) in \(\mathrm{CH}_{r-1}(F(X))\) such that 
  \[
  q_*p^*\alpha_{r-1}=2h^{n-r}_X,~ q_*p^*\beta_{r-1}=5h^{n-r}_X.
  \]
\end{lemma}

\begin{proof}
We first prove the statement by assuming that \(X\) has a \(k\)-point 
\(x\) and a \(k\)-line \(L\) in general position. 

Consider the closed subvariety 
\[
C_x\coloneqq \{[\ell]\in F(X)~|~ x\in \ell\}
\]
of lines meeting with the point \(x\). Since \(x\) is a general point, 
the dimension of \(C_x\) is \(n-3\). Fix a coordinate 
\((x_0\colon, \ldots, \colon x_{n+1})\) of \(\mathbb{P}^{n+1}\). 
By projective transformation we may set \(x=(1\colon 0 \colon ,\ldots, \colon 0)\). 
Then the defining equation of \(X\) is of the form \(x_0^2G_1+x_0G_2+G_3\) 
where \(G_1\) (resp. \(G_2\) and \(G_3\)) is a homogeneous polynomial in 
\(k[x_1,\ldots,x_{n+1}]\) of degree \(1\) (resp. \(2\) and \(3\)). 
Each line in \(X\) meeting with the point \(x\) corresponds
to a point in \(\mathbb{P}(x_1\colon, \ldots, \colon x_{n+1})\) cut out
by the equations \(G_1, G_2\) and \(G_3\). The closed subvariety \(q_*p^*C_x\)
in \(X\) is simply defined by \(G_1=G_2=0\). Hence we have 
\[
q_*p^*[C_x]=2h^2_X\in \mathrm{CH}^2(X).
\] 
Consider a generic linear section 
\(X_{r+2}:=\mathbb{P}^{r+3}\cap X\) of codimension \(n-r-2\) in \(X\)
passing through \(x\). Define
\[
C_{x, X_{r+2}}:=\{[\ell]\in F(X)~|~x\in\ell\subset X_{r+2}\}
\]
the closed subvariety of the lines in \(X_{r+2}\) meeting with 
\(x\in X_{r+2}\). The subavriety \(C_{x, X_{r+2}}\) is thus 
contained in \(F(X_{r+2})\). It follows that
\[
q_*p^*[C_{x, X_{r+2}}]=2i_*h^2_{X_{r+2}}=2h^{n-r}_X
\]
with the inclusion \(i\colon {X_{r+2}}\hookrightarrow X\). Taking 
\(\alpha_{r-1}=[C_{x, X_{r+2}}]\) yields to the first assertion.

For the second assertion, we define \(S_L\subset F(X)\) to be the 
closure of the subset
\[
\{[\ell]\in F(X)~|~\ell\cap L\neq \emptyset, ~\ell\neq L\}
\]
of lines meeting with \(L\). Since \(L\) is in general position,
the dimension of \(S_L\) is \(n-2\). Then \(q_*p^*[S_L]\) is a
divisor class on \(X\). Since \(\dim X\geq 3\) the Picard group
of \(X\) is generated by the hyperplane section class \(h_X\), 
which is the Grothendieck-Lefschetz theorem~\cite[IX, \S 3]{Hartshorne:ample-subvar}.
Hence we can write \(q_*p^*[S_L]=ah_X\) for some \(a\in \mathbb{Z}\). 
To determine \(a\), pick another general line \(L'\) on \(X\). 
The projection formula implies that
\[
a=q_*p^*[S_L]\cdot L'=[S_L]\cdot p_*q^*L'
 =[S_L]\cdot [S_{L'}].
\]
Thus the integer \(a\) indicates the number of lines on \(X\) meeting 
with two lines in general position. It has been proved 
in~\cite[Lemma 3.10]{shen14} or~\cite[\S, 5, Lem. 1.14]{Huy_cubic-2023}
that \(a=5\). Again let \(X_{r+1}\colon=\mathbb{P}^{r+2}\cap X\) be 
a generic linear section of codimension \(n-r-1\) in \(X\) containing
\(L\). Define 
\[
S_{L, X_{r+1}}:=\{[\ell]\in F(X)~|~\ell\subset X_{r+1}, \ell\cap L\neq \emptyset\}.
\]
The subvariety \(S_{L, X_{r+1}}\) is contained in \(F(X_{r+1})\). 
It follows that
\[
q_*p^*[S_{L, X_{r+1}}]=5i_*h_{X_{r+1}}=5h_X^{n-r}, i\colon X_{r+1}\hookrightarrow X.
\]
Hence the second assertion holds by taking \(\beta_{r-1}=[S_{L, X_{r+1}}]\).

Now let \(\sum n_i[L_i]\) be a zero-cycle of degree \(1\) on \(F(X)\) 
by our assumption in the lemma. Each line \(L_i\) is defined over
a field extension \(k\subset k_i\) of finite degree so that 
\(\sum n_i[k_i:k]=1\). Suppose that each \(L_i\) is in general position 
on \(X_{k_i}\). Pick any \(k_i\)-point \(x_i\in L_i\) for each \(i\). 
Through the above argument, there exist cycles
\(\alpha_{r-1, i}, \beta_{r-1,i}\) on \(F(X)_{k_i}\) for each \(i\) 
such that
\[
  q_*p^*\alpha_{r-1,i}=2h^{n-r}_{X_{k_i}},~ q_*p^*\beta_{r-1,i}=5h^{n-r}_{X_{k_i}}.
\]
Let \(\sigma_i: X_{k_i}\to X_k\) and \(\mu_i: F(X)_{k_i}\to F(X)_k\)
be the base change maps. Then the cycle 
\(\alpha_{r-1}:= \sum n_i\cdot {\mu_i}_*\alpha_{r-1,i}\) on 
\(F(X)_k\) maps to 
\[
2h_{X_k}^{n-r}=\sum n_i[k_i:k]\cdot 2h_{X_k}^{n-r}=\sum n_i\sigma_*2h_{X_{k_i}}^{n-r}
\] 
via the cylinder homorphism, and the same argument applies to \(\beta_{r-1, i}\). 

Finally, one can circumvent the genericity assumption of \(L_i\) 
by the moving lemma. It states that for any given non-empty Zariski 
open subset \(U\subset V\) in a smooth \(k\)-variety \(V\), any 
zero-cycle on \(V\) is rationally equivalent to a zero-cycle supported 
in \(U\), see the paper~\cite[Compl\'ement, p. 599]{CT-zero-cycle} by 
Colliot-Th\'el\`ene for the proof. Therefore we can replace the zero-cycle
\(\sum n_i[L_i]\) by a rationally equivalent zero-cycle in general position.
\end{proof}

\begin{theorem}
  \label{surj-cylinder}
Let \(X\subset \mathbb{P}^{n+1}_k\) be a smooth cubic hypersurface of
dimension \(n\geq 3\) over a field \(k\), and let \(F(X)\) be the 
variety of lines on \(X\). Assume that \(X\) contains a one-cycle
of degree \(1\). Then for any \(1\leq r < n-1\), the cylinder map 
\[
q_*p^*\colon \mathrm{CH}_{r-1}(F(X))\to \mathrm{CH}_r(X)
\]
is surjective.
\end{theorem}

\begin{proof}
By Lemma~\ref{lem:equiv-assump}, \(X\) contains a one-cycle of degree 
\(1\) if and only if \(F(X)\) contains a zero-cycle of degree \(1\). 
The assumptions in Corollary~\ref{surj-mod} and Lemma~\ref{hyperplane-class}
are satisfied. Then this theorem immediately follows.
\end{proof}

\begin{remark}
The result dose not hold for divisor classes on cubic
threefolds. It is mentioned in~\cite[\S 5, Rmk., 1.22]{Huy_cubic-2023} 
that for a generic cubic threefold \(X\), the N\'eron-Severi group 
\(\mathrm{NS}(F(X))\) is generated by the primitive class 
\([S_{\ell}]\). Since \(q_*p^*[S_{\ell}]=5h_X\) no divisor on 
\(F(X)\)  sent to \(h_X\) via the cylinder map.
\end{remark}

The following lemma, given by the referee, illustrates the 
equivalence of the conditions for one-cycles on \(X\) and 
zero-cycles on \(F(X)\).

\begin{lemma}\label{lem:equiv-assump}
Let \(X\) be a smooth cubic hypersurface in \(\mathbb{P}^{n+1}\) 
over a field \(k\). Assume \(n \geq 3\). The following are equivalent:
\begin{enumerate}
    \item \(F(X)\) contains a zero-cycle of degree \(1\).
    \item \(X\) contains a one-cycle of degree \(1\).
\end{enumerate}  
\end{lemma}
\begin{proof}
The cylinder homomorphism yields \((1) \Rightarrow (2)\), and 
it remains to show \((2) \Rightarrow (1)\).

When \(k\) is a finite field, by the Lang--Weil estimate, 
\(F(X)\) contains a zero-cycle of degree \(1\).

When \(k\) is infinite, by Bertini’s theorem, there exists 
a smooth linear section \(S \subset X\) of dimension \(2\). 
Over the algebraically closed field \(\bar{k}\), the smooth 
cubic surface \(S_{\bar{k}}\) contains exactly \(27\) lines. 
Note that the degree of a zero-cycle does not change under 
field extensions. Then the Fano scheme \(F(S) \subset F(X)\) 
gives rise to a zero-cycle of degree \(27\) on \(F(X)\).

We now assume that \(X\) contains a one-cycle \(\alpha\) of degree \(1\), 
let \(P\subset X \times F(X)\) be the universal line, and let
\[
\beta := P_{*} \alpha \in \mathrm{CH}_{n-2}(F(X)).
\]
Since \(\dim F(X) = 2n - 4\), we see that $\beta$ has the middle (co)dimension.
We claim that \(\deg(\beta^{2}) = 5\). Again we may assume that \(k\) is algebraically closed. Since \(\mathrm{Pic}(X) = \mathbb{Z}\) is generated by the class of a hyperplane section, we have an isomorphism
\[
\deg : \mathrm{CH}_{1}(X)/\mathrm{num} \xrightarrow{\ \sim\ } \mathbb{Z}.
\]
where \(\mathrm{num}\) is the numerical equivalence. In particular, for any line \(l\) on \(X\), we have \(\alpha = [l] \in \mathrm{CH}_{1}(X)/\mathrm{num}\), which implies
\[
\beta = P_{*} [l] \in \mathrm{CH}_{n-2}(F(X))/\mathrm{num}.
\]
Letting \(l_1, l_2 \subset S\) be any two skew lines, we have
\[
\deg(\beta^{2}) = \deg\big( P_{*}[l_1] \cdot P_{*}[l_2] \big),
\]
and since any line on \(X\) that intersects both \(l_1\) and \(l_2\) 
lies on \(X \cap \mathrm{Span}(l_1, l_2) = S\), the right-hand side 
of the above equation equals the number of lines on \(S\) that 
intersect both \(l_1\) and \(l_2\), which is \(5\). The claim follows.

Since \(5\) and \(27\) are coprime, \(F(X)\) contains a zero-cycle of 
degree \(1\). This completes the proof. 
\end{proof}

\section{Integral Hodge conjectures and Tate conjectures}
\label{sec:IHC-ITC}

\begin{theorem}
  \label{int-Hodge-one-cycl}
  The integral Hodge conjecture holds for one-cycles on the variety
  \(F(X)\) of lines of a smooth complex cubic fourfold \(X\).
\end{theorem}
\begin{proof}
  Suppose that \(\alpha\in \mathrm{H}^6(F(X), \mathbb{Z})\) is an 
  integral Hodge class of type \((3,3)\). Then \(q_*p^*\alpha\in \mathrm{H}^4(X,\mathbb{Z})\) 
  is an integral Hodge class of type \((2,2)\). The integral Hodge 
  conjecture for a smooth cubic \(4\)-fold is proved by Voisin, 
  see~\cite{Voisin:Hodge-conjecture}. Hence there exists a \(2\)-cycle
  \(\gamma\in \mathrm{CH}_2(X)\) such that the cohomology class 
  \([\gamma]=q_*p^*\alpha\). As a result of Theorem~\ref{surj-cylinder}, 
  there exists \(1\)-cycle \(\Gamma\in \mathrm{CH}_1(F(X))\) such that 
  \(q_*p^*\Gamma=\gamma\). By the commutative diagram of the cylinder 
  homomorphisms and cycle class maps
  \[
  \begin{tikzcd}
  \mathrm{CH}_1(F(X))\ar[r, "q_*p^*"] \ar[d] & \mathrm{CH}_2(X) \ar[d] \\
  \mathrm{H}^6(F(X),\mathbb{Z})\ar[r, "q_*p^*"] & \mathrm{H}^4(X,\mathbb{Z}),
  \end{tikzcd}
  \]
  we see \(q_*p^*([\Gamma]-\alpha)=0\). The cylinder map \(q_*p^*\)
  on the cohomology groups is an isomorphism since it is dual to the 
  Abel-Jacobi isomorphism
  \[
  p_*q^*: \mathrm{H}^4(X, \mathbb{Z})\to \mathrm{H}^2(F(X), \mathbb{Z})
  \]
  by Beauville and Donagi~\cite{Beauville-Donagi:cubic4fold}.
  Therefore \([\Gamma]=\alpha\) is an algebraic class.
\end{proof}

Let \(V\) be a smooth proper variety over a field \(k\) which is 
finitely generated over its prime field. Let \(\bar{k}\) be the 
separable closure of \(k\), and \(G_k\) be the absolute 
Galois group \(\mathrm{Gal}({\bar{k}}/k)\). Fix a prime number \(\ell\) 
which is invertible in \(k\). The Tate conjecture states that any
\(G_k\)-invariant class in the \(\ell\)-adic cohomology space
\(\mathrm{H}^{2i}(V_{\bar{k}}, \mathbb{Q}_{\ell}(i))\) is spanned 
by the algebraic classes of codimension \(i\) on \(V_k\).

For an intermediate field \(k\subset k'\subset \bar{k}\) of finite 
degree over \(k\), the Galois group \(G_{k'}: = \mathrm{Gal}(\bar{k}/k')\) 
is an open subgroup of the profinite group \(G_k\). The group \(G_{k'}\) 
acts continuously on the \(\mathbb{Z}_{\ell}\)-module 
\(\mathrm{H}^{2i}_{\et}(V_{\bar{k}}, \mathbb{Z}_{\ell}(i))\). For 
codimension \(i\) algebraic cycles on \(V_{\bar{k}}\), there is
the cycle class map 
\[
\mathrm{cl}^{i}_{\bar{V}} \colon \mathrm{CH}^{i}(V_{\bar{k}})\underset{\mathbb{Z}}{\otimes} \mathbb{Z}_{\ell}
\to \underset{G_{k'}}{\varinjlim}~\mathrm{H}^{2i}_{\et}(V_{\bar{k}}, \mathbb{Z}_{\ell}(i))^{G_{k'}}
\]
where the direct limit is over all intermediate fields 
\(k\subset k'\subset \bar{k}\) of finite degrees over \(k\). Schoen 
proposed the following integral analog of the Tate conjecture for 
one-cycles. 
\begin{conjecture}
The cycle class map \(\mathrm{cl}^i_{\bar{V}}\) is surjective for 
\(i=\dim V-1\).  
\end{conjecture}
The Tate conjecture is true for divisors on the variety \(F(X)\) of 
lines on a smooth cubic fourfold \(X\) over a number field or a finite 
field of characteristic \(p\geq 5\), see~\cite{Shaf-Tate} and~\cite{Tate-K3}. 
Thus the Tate conjecture holds for one-cycles on \(F(X)\) by the hard 
Lefschetz's theorem. Let us prove the integral Tate conjecture for 
one-cycles on \(F(X)\).

\begin{theorem}
  \label{ITC}
  Let \(k\) be a finitely generated field of characteristic different
  from \(2\) and \(3\). Denote by \(F\) the variety \(F(X)\) of lines 
  on a smooth cubic fourfold \(X\) defined over \(k\). Fix a prime 
  \(\ell\) invertible in \(k\), the cycle class map
  \[
  \mathrm{cl}^3_{\bar{F}}\colon \mathrm{CH}_1(F_{\bar{k}})\otimes \mathbb{Z}_{\ell}\rightarrow 
  \underset{U}{\varinjlim}~\mathrm{H}^{6}_{\et}(F_{\bar{k}}, \mathbb{Z}_{\ell}(3))^{U}
  \]
  is surjective, where \(U\) runs over all open subgroups of 
  \(\mathrm{Gal}(\bar{k}/k)\).
\end{theorem}

To imitate the proof of Theorem~\ref{int-Hodge-one-cycl}, 
we need the following lemma.

\begin{lemma}
  \label{etale-Abel-Jacobi}
  Let \(k\) be a number field or a finite field. Denote by 
  \(G_k:= \mathrm{Gal}(\bar{k}/k)\). Fix a prime \(\ell\) 
  invertible in \(k\). Then the cylinder map on the \(\ell\)-adic cohomology 
  \begin{equation}
  \label{eqs:etale-cyl-map}
  \mathrm{H}^6_{\et}(F_{\bar{k}}, \mathbb{Z}_{\ell}(3))
  \rightarrow \mathrm{H}^4_{\et}(X_{\bar{k}}, \mathbb{Z}_{\ell}(2))
  \end{equation}
is a \(G_k\)-equivariant isomorphism.
\end{lemma}
\begin{proof}
Let \(P_{\bar{k}}\) be the \(\mathbb{P}^1\)-bundle over \(F(X)\). 
The projections \(p: P_{\bar{k}}\to F_{\bar{k}}\) and 
\(q: P_{\bar{k}}\to X_{\bar{k}}\) are naturally \(G_k\)-equivariant. 
Then the induced cylinder map~\eqref{eqs:etale-cyl-map} is  
\(G_k\)-equivariant. 

Now we show that~\eqref{eqs:etale-cyl-map} is an isomorphism by cohomology comparison.
Let \(W(k)\) be the ring of Witt vectors of a finite field \(k\). 
It is known any smooth cubic fourfold \(X\) over \(k\) can lift to 
a family \(\mathcal{X}\) of smooth cubic fourfolds over \(W(k)\). 
Denote by \(\mathcal{F}\) the relative variety of lines of 
\(\mathcal{X}\) over \(W(k)\). The fraction field of \(W(k)\) is of
characteristic zero, which can be embedded into the complex numbers 
\(\mathbb{C}\). Hence the base change \(\mathcal{X}_{\mathbb{C}}\) 
is a complex smooth cubic fourfold, and \(\mathcal{F}_{\mathbb{C}}\) 
is the variety of lines on \(\mathcal{X}_{\mathbb{C}}\).

By the smooth and proper base change and comparison theorem, we have
canonical isomomorphisms 
\[
  \mathrm{H}^{2i}_{\et}(X_{\bar{k}}, \mathbb{Z}_\ell(i))\simeq 
  \mathrm{H}^{2i}_{\et}(\mathcal{X}_{\mathbb{C}}, \mathbb{Z}_{\ell}(i))\simeq 
  \mathrm{H}^{2i}(\mathcal{X}_{\mathbb{C}}^{an}, \mathbb{Z}(i))\otimes \mathbb{Z}_\ell.
\]
Cylinder maps commute with the above comparisons since the 
corresponding cylinder maps are defined by the canonical \(\mathbb{P}^1\)-bundles.
Then~\eqref{eqs:etale-cyl-map} is an isomorphism since
the cylinder map on the Betti cohomology
\[
\mathrm{H}^6(\mathcal{X}^{an}_{\mathbb{C}}, \mathbb{Z}(3))\rightarrow 
\mathrm{H}^4(\mathcal{F}^{an}_{\mathbb{C}}, \mathbb{Z}(2))
\]
is an isomorphism.

Suppose that \(k\) is a number field. Then \(\bar{k}\subset\mathbb{C}\)
both are separably closed fields. It follows 
from~\cite[\S VI, Cor. 4.3.]{etale-Milne} that
\[
\mathrm{H}^{2i}_{\et}(X_{\bar{k}}, \mathbb{Z}_{\ell}(i))\simeq 
\mathrm{H}^{2i}_{\et}(X_{\mathbb{C}}, \mathbb{Z}_{\ell}(i)).
\]
Then the assertion for number fields also follows from the cohomology 
comparison.
\end{proof}

The integral Tate conjecture for \(2\)-cycles of cubic fourfolds has
been proved by Charles and Pirutka.
\begin{theorem}\cite[Theorem 1.1]{entiere-Tate}
  \label{entiere-Tate}
  Let \(k\) be finitely generated field of characteristic different 
  from \(2\) and \(3\), and let \(\bar{k}\) be the separable closure 
  of \(k\). Let \(X\) be a smooth cubic fourfold over \(k\). Fix a
  prime \(\ell\) invertible in \(k\), the cycle class map
  \[
  \mathrm{cl}^2_{\bar{X}}\colon \mathrm{CH}^{2}(X_{\bar{k}})\otimes \mathbb{Z}_{\ell}
  \rightarrow \underset{U}{\varinjlim}~\mathrm{H}^{4}_{\et}(X_{\bar{k}}, \mathbb{Z}_{\ell}(2))^{U}
  \]
  is surjective, where the direct limit is over all open subgroups \(U\) 
  of \(\mathrm{Gal}(\bar{k}/k)\).
\end{theorem}

\begin{proof}[Proof of Theorem~\ref{ITC}]
Let \(k_0\) be the prime subfield of \(k\). Let \(R\) be the domain 
of finite type over \(k_0\) such that its fraction field is \(k\). 
After shrinking \(\mathrm{Spec}(R)\) to an open affine subset, we may 
extend the \(X\) defined over \(k\) to a family of smooth cubic fourfolds 
\(\pi\colon \mathcal{X}\to R\). Then the \'etale sheaf 
\(R^4\pi_*\mathbb{Z}_{\ell}(2)\) is locally constant on 
\(\mathrm{Spec}(R)\). Regarding \(X_{\bar{k}}\) as the 
geometric generic fiber of \(\pi\), and \(X_{\bar{k}_0}\)
as the geometric fiber of the central point in \(\mathrm{Spec}(R)\). 
Since \(\pi\) is smooth and proper, the specialization map
\[
\mathrm{H}^4(X_{\bar{k}}, \mathbb{Z}_{\ell}(2))\to \mathrm{H}^4(X_{\bar{k}_0}, \mathbb{Z}_{\ell}(2))
\]
is an isomorphism, see Corollary~\cite[\S VI, Cor. 4.2.]{etale-Milne}. 
The specialization maps for \(X\) and \(F\) commute with the cylinder 
maps because the later are defined by the universal \(\mathbb{P}^1\)-bundle.
Since \(k_0\) is either a finite field or the field of rational numbers,
it implies by Lemma~\ref{etale-Abel-Jacobi} that
\begin{equation}
\label{eqs:cyl-map-fin-gen-field}
\mathrm{H}^{6}_{\et}(F_{\bar{k}}, \mathbb{Z}_{\ell}(3))\to \mathrm{H}^4(X_{\bar{k}}, \mathbb{Z}_{\ell}(2))
\end{equation}
is a \(G_k\)-equivariant isomorphism.

Let \(\alpha\in \mathrm{H}^{6}_{\et}(F_{\bar{k}}, \mathbb{Z}_{\ell}(3))^{U}\)
represents a given cohomology class in the direct limit
\(\underset{U}{\varinjlim} \mathrm{H}^{6}_{\et}(F_{\bar{k}}, \mathbb{Z}_{\ell}(3))^{U}\). 
By Theorem~\ref{entiere-Tate}, there exists a \(2\)-cycle 
\(\Gamma\in \mathrm{CH}^2(X_{\bar{k}})\otimes \mathbb{Z}_{\ell}\)
such that \(\mathrm{cl}^2_{\bar{X}}(\Gamma)=q_*p^*\alpha\) in 
\(\mathrm{H}^4_{\et}(X_{\bar{k}}, \mathbb{Z}_{\ell}(2))^U\). 
By the surjectivity of the cylinder map, there exists a one-cycle
\(\gamma\in \mathrm{CH}_1(F_{\bar{k}})\otimes \mathbb{Z}_{\ell}\) 
such that \(q_*p^*\gamma=\Gamma\). Since \eqref{eqs:cyl-map-fin-gen-field}
is an isomorphism we have \(\mathrm{cl}^3_{\bar{F}}(\gamma)=\alpha\).
\end{proof}

\end{document}